\documentclass{article}
\usepackage{amsmath}
\usepackage{amsfonts}
\usepackage{graphicx}
\usepackage[style=numeric,natbib=true,backend=bibtex]{biblatex}
\bibliography{analysis}
\include{minilib}
\include{subsup}
\newcommand{\defmeso}[2]{\newsubsupcommand{#1}{#2}[][\varepsilon]}
\newcommand{\defmesoi}[2]{\newsubsupcommand{#1}{#2}[i][\varepsilon]}
\defineDomain{Time}{t}{(0,T)}
\renewInt{Time}{\int\limits_0^T}
\defineDomain{A}{x}{\Omega^\varepsilon}
\defineBoundarySegment{G}{A}{Grain}{\Gamma^\varepsilon}
\defineBoundarySegment{ANeumann}{A}{Neumann}{\Gamma_N^\varepsilon}
\defineBoundarySegment{ARobin}{A}{Robin}{\Gamma_R^\varepsilon}
\defineBoundarySegment{ATRobin}{A}{TRobin}{\Sigma_R}
\defineBoundarySegment{ATDirichlet}{A}{TDirichlet}{\Sigma_D}
\defineCrossProduct{TA}{Time}{A}
\defineCrossProduct{TG}{Time}{G}
\newcommand\tensorIso{0.75 & 0.171476\\0.171476 & 0.75}
\newcommand\tensorAniso{0.817467 & 0.0786338\\0.214942 & 0.817467}
\newcommand\tensorIsotor{1.0000 &  0.2286\\0.2286 & 1.0000}
\newcommand\tensorAnisotor{ 0.9617 &  0.0925\\  0.2529 & 0.9617}
\newcommand\tempBase{\theta}
\newcommand\collBase{u}
\newcommand\depBase{v}
\defmeso{\temper}{\tempBase}
\defmeso{\coller}{\collBase}
\defmeso{\deper}{\depBase}
\defineUnknown{temp}{\temper}{TA}
\defineUnknown{coll}{\coller}{TA}
\defineUnknown{dep}{\deper}{TG}
\inheritUnknown{coll}{coll_i}{\coller_i}
\inheritUnknown{dep}{dep_i}{\deper_i}
\defmesoi{\colli}{u}
\defmesoi{\depi}{v}

\newcommand{\bigO}[1]{\mathcal{O}(\varepsilon^#1)}
\defmeso{\cond}{\kappa}
\defmeso{\diff}{d}
\attachFlux{temp}{\cond}
\attachFlux{coll}{\diff}
\attachFlux{coll_i}{\diff_i}

\defmeso{\dca}{a}
\defmeso{\dcb}{b}

\newcommand{\Rdim}[1]{\mathbb{R}^#1}
\newcommand{\Zdim}[1]{\mathbb{Z}^#1}
\newcommand{\Runit}[1]{\vec{e}_#1}

\newcommandx\newcommandWithIndex[3]{
  \constructor{#1}{#2#3}
  \set{#1}{body}{#2}
  \set{#1}{index}{#3}
}
\makeatletter
\newcommand{\closure}[1]{
  \@ifundefined{#1@body}
  {\overline{#1}}
  {\overline{\@nameuse{#1@body}}\@nameuse{#1@index}}
}
\makeatother
\defineDomain{Domain}{x}{\Omega}
\defineBoundarySegment{DTNeumann}{Domain}{TNeumann}{\Gamma_N^\theta}
\defineBoundarySegment{DTRobin}{Domain}{TRobin}{\Gamma_R^\theta}
\defineBoundarySegment{DTDirichlet}{Domain}{TDirichlet}{\Gamma_D^\theta}
\defineBoundarySegment{DUNeumann}{Domain}{UNeumann}{\Gamma_N^u}
\defineBoundarySegment{DURobin}{Domain}{URobin}{\Gamma_R^u}
\defineBoundarySegment{DUDirichlet}{Domain}{UDirichlet}{\Gamma_D^u}
\newcommand{\cellUnit}{Y}
\newcommandWithIndex{cellFluid}{\cellUnit}{_1}
\newcommandWithIndex{cellSolid}{\cellUnit}{_0}
\newcommand{\cellSolidBoundary}{\Gamma}
\newcommandWithIndex{volSkeleton}{\Domain}{_0^\varepsilon}
\newcommandWithIndex{volPore}{\Domain}{^\varepsilon}
\newcommand{\surfSkeleton}{\Gamma^\varepsilon}

\newcommandx{\inRange}[2][2=]{\ifstrequal{#2}{}
  {\in\{1,\ldots,#1\}}
  {\in\{#1,\ldots,#2\}}
}
\newcommand\collu{u}
\newcommand\depu{v}

\newcommand\Dca{A}
\newcommand\Dcb{B}

\newcommand\intG[1]{\int_{\Gamma}#1\,d\sigma(y)}

\newcommand\Thiele{\Lambda}
\newcommand\Biot{Bi}
\newcommand\porosity{\phi}

\begin{document}
\title{Multiscale Modeling of Colloidal Dynamics in Porous Media: Capturing Aggregation and Deposition Effects}


\date{ }
\maketitle
\centerline{\scshape Oleh Krehel}
\medskip
{\footnotesize
  \centerline{Department of Mathematics and Computer Science}
    \centerline{CASA - Center for Analysis, Scientific computing and Engineering}
  \centerline{Eindhoven University of Technology}
  \centerline{5600 MB, PO Box 513, Eindhoven, The Netherlands}
 } 

\medskip

\centerline{\scshape Adrian Muntean}
\medskip

{\footnotesize
  \centerline{Department of Mathematics and Computer Science}
  \centerline{CASA - Center for Analysis, Scientific computing and Engineering}
    \centerline{ICMS - Institute for Complex Molecular Systems}
  \centerline{Eindhoven University of Technology}
  \centerline{5600 MB, PO Box 513, Eindhoven The Netherlands}
}

\medskip
\centerline{\scshape Peter Knabner}
\medskip

{\footnotesize
  \centerline{Department of Mathematics}
  \centerline{Friedrich-Alexander University of Erlangen-Nuremberg}
  \centerline{Cauerstr. 11, Erlangen 91058, Germany}
}

\bigskip

\begin{abstract}
  We investigate the influence of multiscale aggregation and
  deposition on the colloidal dynamics in a saturated porous medium.
  At the pore scale, the aggregation of colloids is modeled by the
  Smoluchowski equation. Essentially, the colloidal mass is
  distributed between different size clusters. We treat these clusters
  as different species involved in a diffusion-advection-reaction
  mechanism. This modeling procedure allows for different material
  properties to be varied between the different species, specifically
  the rates of diffusion, aggregation, deposition as well as the
  advection velocities. We apply the periodic homogenization procedure
  to give insight into the effective coefficients of the upscaled model
  equations.  Benefiting from direct access to microstructural
  information, we capture by means of 2D numerical simulations the
  effect of aggregation on the deposition rates recovering this way
  both the blocking and ripening regimes reported in the literature.

\end{abstract}

\section{Introduction}
\label{sec:introduction}

Colloids are particles with size ranging approximately from $1$ to
$1000$ nm in at least one dimension. They play a significant
functional role in a number of technological and biological
applications, such as waste water treatment, food industry, printing,
design of drug delivery; see
e.g. \cite{Nordbotten,rosenholm2010towards}.  The existing literature
on colloids and their dynamics is huge. Here we only mention that the
self-assembly of collagen structures (basic component of the mechanics
of the human body) together with secondary nucleation effects have
recently been treated in \cite{Lith2013}, starting off from an
interacting particle system for colloids. A detailed discussion of the main principles of
aggregation mechanisms can be found in \cite{peukert2005control}, while a
thorough analysis of the aggregation in terms of ordinary differential
equations can be found e.g. in \cite{camejo2012existence}.

The central topic of this paper is the treatment of the aggregation of
colloids in porous media (particularly, soils) that has been recently
shown to be a dominant factor in estimating contaminant transport; see
\cite{totsche2004mobile}. Essentially, one supposes that that the
presence of colloidal aggregation strongly affects the deposition
rates on the pore (grain) boundary.  Similar aggregation (group
formation, cooperation) patterns can emerge also in pedestrian flows
strongly affecting their viscosity \cite{muntean2014pedestrians}.
Previous investigations on contaminant dynamics in soils, yet not
accounting explicitly for aggregation, can be found, for instance, in
\cite{knabner1996modeling} and \cite{totsche1996modeling}.

Our aim here is to study the influence of multiscale aggregation and
deposition on the colloidal dynamics in a saturated porous medium
mimicking a column experiment performed by Johnson, Sun and Elimelech
and reported in \cite{johnson1996colloid}. For more information on
this experimental context, we refer the reader also to
Refs. \cite{liu1995colloid,johnson1995dynamics}. To get more theoretical insight in this column experiment, we proceed as follows:
As departure point, we assume that at the pore scale we can model the
aggregation of colloids by the Smoluchowski equation. Consequently,
the colloidal mass is distributed between different size clusters. We
treat these clusters as different species involved in a coupled
diffusion-advection-reaction system. This modeling procedure allows
for different material properties to be varied between the different
species, specifically the rates of diffusion, aggregation, deposition
as well as the advection velocities. As next step, we apply the
periodic homogenization methodology to give insight into the effective
coefficients of the upscaled model equations. Finally, for a set of
reference parameters, we solve the upscaled equations for different
choices of microstructures and investigate the influence of
aggregation on both transport and deposition of the colloidal mass,
validating in the same time our methodology and numerical platform by
means of the results from \cite{johnson1996colloid}.

The outline of the paper is as follows: In
Section~\ref{sec:microscopic-model} we set up a microscopic pore-scale
model for aggregation, diffusion and deposition of populations of colloidal
particles.  In Section~\ref{sec:nondim} the microscopic model is
nondimensionalized.  One of the small dimensionless numbers pointed out therein
(denoted by $\varepsilon$) connects a ratio of characteristic time scales of the process to the relevant microscopic and macroscopic
length scales arising in the system.  In
Section~\ref{sec:macroscopic-model} we use the concept of two-scale
asymptotic expansions to obtain in the limit of small $\varepsilon$ an
equivalent macroscopic model together with the corresponding effective
coefficients.  We conclude the paper with a few numerical multiscale
experiments and discussions on further work
(cf. Section~\ref{sec:numerics} and Section \ref{Discussion}).

\section{Microscopic model}
\label{sec:microscopic-model}

The foundations of the modeling of colloids aggregation and fragmentation were laid down in the
classical work of Smoluchowski \cite{smoluchowski1917versuch}.  A nice
overview can be found, for instance, in
\cite{elimelech1995particle}. The role of this section is to introduce
our modeling Ansatz on the second order kinetics describing the
colloidal cluster growth and decline, the functional structure of the
deposition rate, as well as the assumptions on the microscopic
diffusion coefficients for the clusters.

\subsection{Aggregation and fragmentation of clusters}\label{aggregation}

We assume that the colloidal population consists of identical
particles, called primary particles, some of which form aggregate
particles that are characterized by the number of primary particles
that they contain -- i.e. we have $u_1$ particles of size $1$, $u_2$
particles of size $2$, etc. We refer to each particle of size $i$ as a
member of the $i^{th}$ species (or of the $i-$cluster).

The fundamental assumption behind this modeling strategy is that
aggregation can be perceived as a second-order rate process, i.e. the rate of collision
is proportional to concentrations of the colliding species. Thus
$A_{ij}$ -- the number of aggregates of size $i + j$ formed from the
collision of particles of sizes $i$ and $j$ per unit time and volume,
equals:
\begin{align}
  \label{eq:second-order-rate}
  A_{ij} &:= \gamma_{ij} u_i u_j\text{, with}\\
  \label{eq:second-order-rate-coeff}
  \gamma_{ij} &:= \alpha_{ij} \beta_{ij}.
\end{align}
Here $\beta_{ij}$ is the collision kernel -- rate constant determined
by the transport mechanisms that bring the particles in close contact,
while $\alpha_{ij} \in [0,1]$ is the collision efficiency -- the fraction of
collisions that finally form an aggregate. The coefficients
$\alpha_{ij}$ are determined by a combination of particle-particle
interaction forces, both DLVO (i.e. double-layer repulsion and van der
Waals attraction) and non-DLVO, e.g. steric interaction forces (see
\cite{derjaguin1941theory}, \cite{hamaker1937london}).

A typical choice for $\alpha_{ij}$ and $\beta_{ij}$ can be found in for
instance in \cite{krehel2012flocculation}.  The interaction rates
(written in the spirit of balance of populations balances as reaction
rates) should then satisfy
\begin{equation}
  \label{eq:PBE}
  R_i(u)=\frac{1}{2}\sum_{i+j=k}\alpha_{ij}\beta_{ij}u_iu_j
  -u_k\sum_{i=1}^{\infty}\alpha_{ki}\beta_{ki}u_i,
\end{equation}
where $u=(u_1,\ldots,u_N,\ldots)$ is the vector of the concentrations
for each size class $i\in\{1,\dots, N\}$ for a fixed choice of $N$.

\subsection{Diffusion coefficients for clusters}\label{diffusion}

We take the diffusivity $d_1$ of the monomers  as a baseline.
All the other diffusivities are here assumed to depend on $d_1$ in agreement with the
Einstein-Stokes relation
\begin{equation}
  \label{eq:einstein-stokes}
  d_i=\frac{kT}{6\pi\eta r_i}.
\end{equation}
The cluster diffusion coefficients $d_i$ arising in
(\ref{eq:einstein-stokes}) are designed for the diffusion of spherical
particles through liquids at low Reynolds number. In
(\ref{eq:einstein-stokes}), $T$ denotes the absolute temperature, $k$
is the Boltzmann factor, $\eta$ is the dynamic viscosity, while $r_i$
is the aggregate ($i$-mer, $i$-cluster) radius.  Note the following
dependence of the aggregate radius $r_i$ on the number of monomers
contained in the $i$-cluster:
\begin{equation}
  \label{eq:aggregate-radius}
  r_i=i^\frac{1}{D_F}r_1,
\end{equation}
with $D_F$ being a dimensionless parameter called the fractal
dimension of the aggregate \cite{meakin1987fractal}.  $D_F$ indicates
how porous the aggregate is. For instance, a completely non-porous
aggregate in three dimensions, such as coalesced liquid drops, would have $D_F=3$.  Combining
(\ref{eq:einstein-stokes}) and (\ref{eq:aggregate-radius}), we obtain:
\begin{equation}
  \label{eq:diffusion-i}
  d_i=\frac{1}{i^\frac{1}{D_F}}d_1.
\end{equation}

\subsection{Deposition rate of colloids on grain surfaces}\label{deposition}

The colloidal species $u_i$, defined in $\Omega$ (see
Figure~\ref{fig:porous}), can deposit on the grain boundary of the solid
matrix $\Gamma\subset\partial \Omega$, transforming into an immobile species
$v_i$, defined on $\Gamma$. This means that the colloids of different
size can be present both in the bulk and on the boundary.  The
boundary condition for $\Gamma$ then looks like:
\begin{equation}
  \label{eq:deposition-bc}
  -d_i\nabla u_i\cdot n=F_i(u_i,v_i).
\end{equation}
At this stage, we assume the deposition rate $F_i$ to be linear, namely we take
\begin{equation}
  \label{eq:deposition-bc-linear}
  F_i(u_i,v_i)=a_iu_i - b_iv_i,
\end{equation}
this resembles the structure of Henry's law acting in the context of
gas exchange at liquid interfaces \cite{battino1966solubility}.

\subsection{Setting of the microscopic model equations}

Collecting the modeling assumptions from Section \ref{aggregation},
Section \ref{diffusion}, and Section \ref{deposition}, we see that the
microscopic system to be tackled in this context is as follows:

Find $(u_1,\dots, u_N, v_1,\dots,v_N$) satisfying
\begin{align}
  \label{eq:micro-start}
  & \partial _t u_i + \nabla \cdot (-d_i\nabla u_i) = R_i(u) && \text{ in }\Omega,\\
  & \partial _t v_i = a_iu_i - b_iv_i && \text{ on }\Gamma,\\
  \intertext{with the boundary conditions}
  & - d_i\nabla u_i\cdot n = a_iu_i - b_iv_i && \text{ on }\Gamma,\\
  & - d_i\nabla u_i\cdot n = 0              && \text{ on }\Gamma_N,\\
  &   u_i = u_{iD}        && \text{ on }\Gamma_D,\\
  \intertext{and the initial conditions}
  & u_i(0,x)=u_i^0(x) && \text{ for } x\in \Omega,\\
  \label{eq:micro-end}
  & v_i(0,x)=v_i^0(x) && \text{ for } x\in\Gamma.
\end{align}

\section{Nondimensionalization}
\label{sec:nondim}

Let $\tau$, $\chi$, $d$, $\collu_0$, $\depu_0$, and $a_0$ be reference
quantities. We choose the scaling $t:=\tau\tilde{t}$,
$x:=\chi\tilde{x}$, $d_i:=d\tilde{d}_i$,
$\collu_i:=\collu_0\tilde{\collu}_i$,
$\depu_i:=\depu_0\tilde{\depu}_i$, $a_i:=a_0\tilde{a}_i$, and
$b_i:=\frac{a_0u_0}{v_0}\tilde{b_i}$.  As reference quantities, we select $\chi:=L$,
$d:=d_1$, $\collu_0:=\max\{u_{i0}, u_{iD}:i\in\{1,\dots,N\}\}$, and
$\depu_0:=\max\{v_{i0}:i\in\{1,\dots,N\}\}$.

Note that we need to distinguish between $\collu_0$ and $\depu_0$
since they have different dimensions, i.e. volume and surface
concentration, respectively.  After substituting these scaling relations into
(\ref{eq:micro-start})-(\ref{eq:micro-end}) and  dropping the tildes, we obtain:
\begin{align}
  \partial _tu_i+\frac{\tau d}{L^2}\nabla \cdot (-d_i\nabla u_i)=\tau u_0R_i(u)\\
  -d_i\nabla u_i\cdot n=\frac{a_0L}{d}(a_iu_i-b_iv_i)\\
  \partial _tv_i=\frac{\tau a_0}{v_0}u_0(a_iu_i-b_iv_i).
\end{align}
This nondimensionalization procedure involves three relevant dimensionless numbers.
We denote by $\varepsilon$ our first dimensionless number, viz.
\begin{equation}\label{epsilon_number}
\varepsilon:=\frac{a_0L}{d}.
\end{equation}
For our particular scenario, the dimensionless number $\varepsilon $ takes a small value (here
$\varepsilon\approx 7.61e-7$).  We will relate it  in Section \ref{sec:macroscopic-model} to a ratio of characteristic
micro-macro length scales. We
refer to $\varepsilon $ as the {\em homogenization parameter}. Furthermore, we choose to scale
the time variable in the system by the characteristic time scale of
diffusion $\tau:=\frac{L^2}{d}$ of the fastest species (i.e. the monomers).  This particular choice of time scale
leads to two further dimensionless numbers:
\begin{itemize}
\item the {\em Thiele modulus}
\begin{equation}\label{Thiele_number}
\Thiele:=\frac{L^2}{d}u_0
\end{equation}
\item the {\em Biot number}
\begin{equation}\label{Biot_number}
\Biot:=a_0\frac{L^2}{d}\frac{u_0}{v_0}.
\end{equation}
\end{itemize}
According to our reference parameters, we estimate that $\Thiele=3.8934e21$ and $\Biot=7.6914e-08$. The order of magnitude of the Thiele modulus $\Thiele$ indicates that the characteristic reaction time is very small compared to the characteristic time of monomers diffusion, the overall reaction-diffusion process being with this scaling in its fast reaction regime. The order of magnitude of the Biot number $\Biot$ points out the slow deposition regime. Essentially, since $\frac{Lu_0}{v_0}=\mathcal{O}(1)$, we have $\Biot=\mathcal{O}(\varepsilon)$. To remove a proportionality constant in the scaled boundary condition (\ref{bibi}),  we take $L:=\frac{v_0}{u_0}$.

Finally, we obtain the following dimensionless system of governing equations:

\begin{align}
  \label{eq:micro-nodim-start}
  & \partial _t u_i + \nabla \cdot (-d_i\nabla u_i) = \Thiele R_i(u) && \text{ in }\Omega,\\
  & \partial _t v_i = \Biot(a_iu_i-b_iv_i)             && \text{ on }\Gamma,\\
  \intertext{with the boundary conditions}
  & - d_i\nabla u_i\cdot n = \varepsilon (a_iu_i-b_iv_i)          && \text{ on }\Gamma,\label{bibi}\\
  & - d_i\nabla u_i\cdot n = 0                        && \text{ on }\Gamma_N,\\
  &   u_i(t,x) = \frac{u_D(t,x)}{u_0}      && \text{ on }\Gamma_D,\\
  \intertext{and the initial conditions}
  & u_i(0,x)=\frac{u_i^0(x)}{u_0}          && \text{ for } x\in\Omega,\\
  \label{eq:micro-nodim-end}
  & v_i(0,x)=\frac{v_i^0(x)}{v_0}          && \text{ for } x\in\Gamma.
\end{align}


\section{Derivation of the macroscopic model}
\label{sec:macroscopic-model}

In this section, we suppose that our porous medium has an internal
structure that can be sufficiently well approximated by an array of
periodically-distributed microstructures. For this situation, starting
off from a partly dissipative model for the dynamics of large
populations of interacting colloids at the pore level (i.e. within the
microstructure), we derive upscaled equations governing the
approximate macroscopically observable behavior. To do this, we employ
the technique of periodic homogenization; see, for instance,
\cite{BLP,Marchenko,CPS}.  In what follows, we apply the technique in
an algorithmic way, giving complete and explicit
calculations. 

\subsection{Colloid dynamics in structured media. The periodic homogenization procedure}

\begin{center}
    \begin{figure}[!ht]
      \centering
      \begin{tabular}{ll}
        \includegraphics[width=0.45\textwidth]{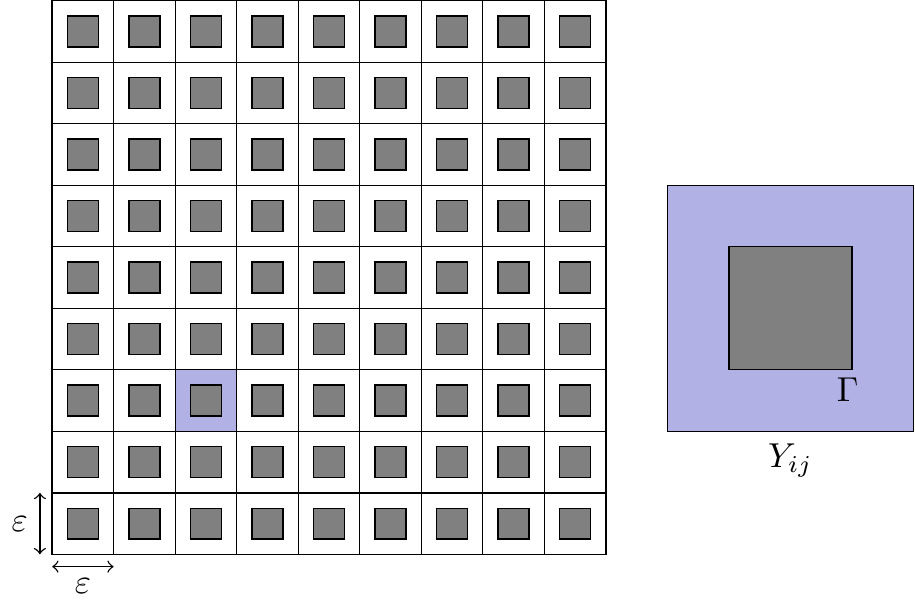}
        &
          \includegraphics[width=0.45\textwidth]{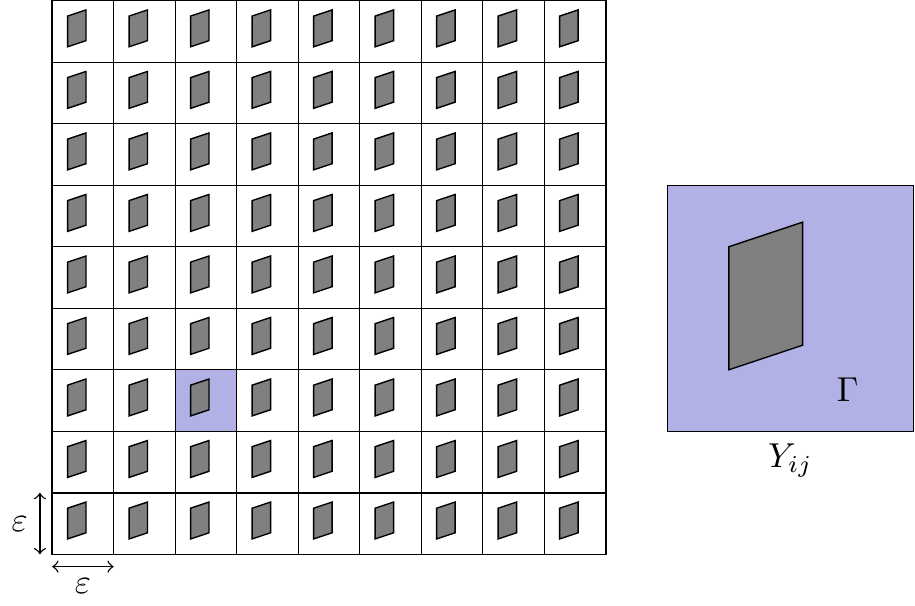}
      \end{tabular}
      \caption{Microstructure of $\Omega^\varepsilon$. Left: isotropic case; Right:  anisotropic case. Here  $Y_{ij}$ is the
        periodic cell.}
      \label{fig:porous}
    \end{figure}
\end{center}

The porous medium $\Omega^\varepsilon$ that we consider  is modeled here as
a composite periodic structure with $\varepsilon>0$ as a small scale
parameter, which relates the the pore length scale to the domain
length scale.
$\Omega^\varepsilon$ is depicted in Figure~\ref{fig:porous}. We assume in this context that this scale parameter is of the same order of magnitude as $\varepsilon$ introduced in (\ref{epsilon_number}).
Note in Figure~\ref{fig:porous} the periodic array of cells approximating the porous media under consideration. Each element
is a rescaled  (by $\varepsilon$) and translated copy of the standard cell $Y$.

\begin{table}[h!]
  \begin{tabular}{ll}
  \hline
  \\
    $\Time$        & $=$ time interval of interest\\
    $\Domain$      & $=$ bounded domain in $\Rdim{n}$\\
    $\partial \Domain$     & $=\Gamma_R\cup \Gamma_N$
                     piecewise smooth boundary of $\Domain$, $\Gamma_R\cap \Gamma_N=\emptyset$\\
    $\Runit{i}$    & $=$ $i$th unit vector in $\Rdim{n}$ ($n=2$ or $n=3$)\\
    $\cellUnit$    & $=\{\sum_{i=1}^n\lambda_i\Runit{i}:\: 0<\lambda_i<1\}$ unit cell in $\Rdim{n}$\\
    $\cellSolid$   & $=$ open subset of $\cellUnit$ that represents the solid grain\\
    $\cellFluid$   & $=\cellUnit \setminus  \closure{cellSolid}$ \\
    $\cellSolidBoundary$    & $=\partial \cellSolid$ piecewise smooth boundary of $\cellSolid$\\
    $X^k$          & $=X+\sum_{i=1}^nk_i\Runit{i}$, where $k\in \Zdim{n}$ and $X\subset\cellUnit$\\
    \\
    \hline
  \end{tabular}
  \caption{$\varepsilon$-independent objects.}
\end{table}


\begin{table}[!ht]
  \begin{tabular}{ll}
  \hline\\
    $\volSkeleton$     & $=\cup \{\varepsilon \cellSolid^k:\: \cellSolid^k\subset \A, k\in \Zdim{n}\}$ array of pores\\
    $\volPore$         & $=\Domain\setminus\closure{volSkeleton}$ matrix skeleton
    \\
    $\surfSkeleton$ & $=\partial \volSkeleton$ pore boundaries\\
    \\
    \hline
  \end{tabular}
  \caption{$\varepsilon$-dependent objects.}
\end{table}
As customary in periodic homogenization applications, we introduce the fast variable $y:=x/\varepsilon$ and let all the unknowns be represented
by the following expansions:
\newcommand\asymptoticExpansion[1] {
  #1(x)&:=#1_0(x,y)+\varepsilon#1_1(x,y)+\varepsilon^2#1_2(x,y)+\bigO{3}
}
\begin{equation}
  \label{eq:asymptotic-expansion}
  \begin{cases}
    \asymptoticExpansion{\coll},\\
    \asymptoticExpansion{\dep}.
  \end{cases}
\end{equation}
\undef\asymptoticExpansion
The asymptotic expansions (\ref{eq:asymptotic-expansion})
can be justified by means of the concept of two-scale convergence by
Nguetseng and Allaire; see Ref. \cite{Oleh_analysis} for the
mathematical analysis of a more complex case including also thermal
effects, and \cite{hornung1991diffusion} for a closely related
scenario. 

Now, taking into account the chain rule
$\nabla :=\nabla _x+\frac{1}{\varepsilon}\nabla _y$, we get:
\newcommand\asymptoticExpansion[2] {
  #1#2&=
  \varepsilon^{-1}#1_y#2_0+
  \varepsilon^{0}(#1_x#2_0+#1_y#2_1)+
  \varepsilon^{1}(#1_x#2_1+#1_y#2_2)+
  \bigO{2}.
}

\begin{align*}
  \asymptoticExpansion{\nabla }{\colli}\\
  \asymptoticExpansion{\nabla }{\depi}
\end{align*}
\undef\asymptoticExpansion
This gives us the following diffusion term:
\begin{align*}
  \nabla \cdot (\diff_i(y)\nabla \colli) & = \varepsilon ^{-2}\nabla _y\cdot (\diff_i(y)\nabla _y\colli_0)\\
                        & +\varepsilon ^{-1}
                          (
                          \diff_i(y)\nabla _x\cdot \nabla _y\colli_0
                          +\nabla _y\cdot (\diff_i(y)\nabla _x\colli_0)
                          +\nabla _y\cdot (\diff_i(y)\nabla _y\colli_1)
                          )\\
                        & +\varepsilon ^{0}
                          (
                          \diff_i(y)\Delta \colli_0
                          +\diff_i(y)\nabla _x\cdot \nabla _y\colli_1\\
                         &\qquad +\nabla _y\cdot (\diff_i(y)\nabla _x\colli_1)
                          +\nabla _y\cdot (\diff_i(y)\nabla _y\colli_2)
                          )+\bigO{1}.
\end{align*}
Collecting the terms with $\varepsilon ^{-2}$ gives:
\begin{align*}
  \nabla _y\cdot (\diff_i(y)\nabla \colli_0)=0.
\end{align*}
Recalling that this PDE with periodic boundary conditions has a solution
unique up to a constant, we get $\colli_0=\colli_0(x)$.
Consequently, we have $\nabla _y\colli_0=0$.

The terms with $\varepsilon ^{-1}$ can be arranged as
\begin{equation}
  \label{eq:eps1}
  \nabla _y\cdot (\diff_i(y)\nabla _y\colli_1)=-\nabla _y\diff_i(y)\cdot \nabla _x\colli_0.
\end{equation}

Let $w_j(y)$ solve the following {\em cell problem} endowed with periodic boundary conditions:
\begin{align}
  \label{eq:cell-problem}
  \nabla _y\cdot (\diff_i(y)\nabla w_j)=-(\nabla \diff_i(y))_j&&j\inRange{d}, y\in Y
\end{align}
Using (\ref{eq:cell-problem}), we can express the first order term in
(\ref{eq:asymptotic-expansion}) as:
\begin{equation}
  \label{eq:eps2}
  \colli_1(x,y)=w(y)\cdot \nabla \colli_0(x)+\colli_1(x),
\end{equation}
where the function $\colli_1(x)$ does not depend on the variable $y$. Note that
\begin{equation}
  \label{eq:eps3}
  \nabla _y\colli_1(x,y)=\nabla w(y)\cdot \nabla \colli_0(x).
\end{equation}

The terms with $\varepsilon ^{0}$ give:
\begin{align*}
  \partial _t\colli_0=&\diff_i(y)\Delta \colli_0+\diff_i(y)\nabla w(y):\nabla \nabla \colli_0\\
              &+\nabla _y\cdot (\diff_i(y)\nabla _x\colli_1+\diff_i(y)\nabla _y\colli_2) + \Thiele R_i(\coll_0).
\end{align*}
Integrating over $Y$ and noting that $|Y|=1$ yield:
\begin{align}
  \label{eq:eps4}
  \partial _t\colli_0=\bar{\mathbb{D}}_i:\nabla \nabla \colli_0-\int_{\partial Y}\diff_i(y)(\nabla _x\colli_1 + \nabla _y\colli_2)\cdot n d\sigma(y) + \Thiele R_i(\coll_0).
\end{align}
The upscaled diffusion tensors $\bar{\mathbb{D}}_i:=[\bar{D}_{ijk}]$ reads:
\begin{align}
  \label{eq:upscaled-diffusion-tensor}
  &\bar{D}_{ijk}=\int_{Y}d_i(y)(\delta_{jk}+\nabla _yw_i(y))dy&&i\inRange{N};j,k\inRange{d}.
\end{align}
Because of the periodic boundary conditions, the active part of $\partial Y$ is only $\Gamma$.
Here we have:
\begin{align}
  \label{eq:eps5}
  \partial _t\colli_0 =
  \bar{\mathbb{D}}_i:\nabla \nabla \colli_0 -
  \int_{\Gamma}\diff_i(y)(\nabla _x\colli_1 + \nabla _y\colli_2)\cdot n d\sigma(y) +
  \Thiele R_i(\coll_0).
\end{align}
The boundary term in (\ref{eq:eps5}) can be expressed recalling the
corresponding deposition boundary condition:
\begin{equation}
  \label{eq:deposition-bc}
  -\diff_i\nabla \coll_i\cdot n=\varepsilon (\dca_i\coll_i-\dcb_i\dep_i)
\end{equation}
Using the prescribed asymptotic expansions, (\ref{eq:deposition-bc}) becomes:
\begin{align*}
  -&\diff_i(y)(\varepsilon^{-1}\nabla _y\colli_0+\varepsilon^0(\nabla _x\colli_0+\nabla _y\colli_1)
  +\varepsilon^1(\nabla _x\colli_1+\nabla _y\colli_2))\cdot n\\
  &=\dca_i(y)(\varepsilon^1\colli_0+\varepsilon^2\colli_1)
  -\dcb_i(y)(\varepsilon^1\dep_0+\varepsilon^2\dep_1)+\bigO{2}.
\end{align*}
Consequently, we obtain
\begin{align*}
  -\diff_i(y)(\nabla _x\colli_1+\nabla _y\colli_2)\cdot n=\dca_i(y)\colli_0-\dcb_i(y)\dep_0.
\end{align*}
Finally, the upscaled equation for $\colli$ reads:
\begin{equation}
  \label{eq:upscaled-coll-1}
  \partial _t\collu_i-\nabla \cdot (\bar{\mathbb{D}}_i\nabla \collu_i)+\Dca_i\collu_i-\Dcb_i\depu_i=\Thiele R_i(\mathbf{\collu}).
\end{equation}
Note that the microscopic surface exchange term turns as
$\varepsilon\to 0$ into the macroscopic bulk term
$\Dca_i\collu_i-\Dcb_i\depu_i$.  Furthermore, the upscaled equation
for $\depi$ is
\begin{equation}
  \partial _t\depu_i=\Dca_i\collu_i-\Dcb_i\depu_i,
\end{equation}
where the effective constants $\Dca_i$ and $\Dcb_i$ are defined by
\begin{equation}
\Dca_i:=\Biot\intG{a_i(y)}
\end{equation}
and
\begin{equation}
  \Dcb_i:=\Biot\intG{b_i(y)}.
\end{equation}

Summarizing, the upscaled system describing the macroscopic dynamics
of the colloids is:
\begin{align}
  \label{eq:macro.model-beg}
  &\partial _t\collu_i-\nabla \cdot (\bar{\mathbb{D}}_i\nabla \collu_i)+\Dca_i\collu_i-\Dcb_i\depu_i=\Thiele R_i(\mathbf{\collu})&&\text{ in }\Omega,i\inRange{N}\\
  &\partial _t\depu_i=\Dca_i\collu_i-\Dcb_i\depu_i&&\text{ in }\Omega,i\inRange{N}\\
  &d_i\nabla \collu_i=f_i&&\text{ on }\Gamma_R,i\inRange{N}\\
  &\collu_i=u_{iD}&&\text{ on }\Gamma_D,i\inRange{N}\\
  &\collu_i(\cdot ,0)=\collu_i^0&&\text{ in }\Omega,i\inRange{N}\\
  \label{eq:macro.model-end}
  &\depu_i(\cdot ,0)=\depu_i^0&&\text{ in }\Omega,i\inRange{N}.
\end{align}


\subsection{Computation of the effective diffusion tensors $\bar{\mathbb{D}}_i=\bar{D}_{ijk}$}

We rely on equation (\ref{eq:upscaled-diffusion-tensor}) to
approximate the main effective transport coefficients -- the effective
diffusion tensors $\bar{D}_{ijk}$ responsible in this scenario for the
transport of the $N$ species of colloids. See
Table~\ref{tbl:upscaled-diffusion} for a calculation example (notice
the symmetry of the tensors corresponding to the isotropic case).

\begin{figure}[h!]
  \centering
  \includegraphics[width=0.8\textwidth]{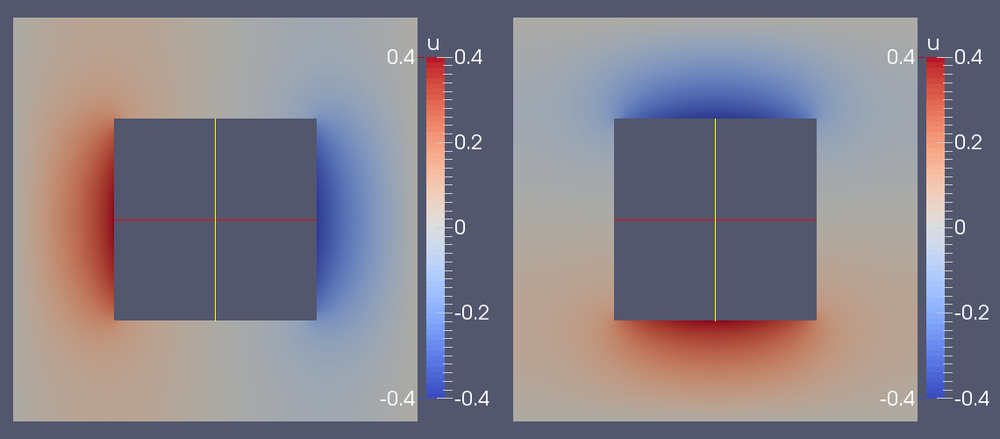}
  \caption{Solutions to the cell problems that correspond to isotropic periodic geometry
    (Figure~\ref{fig:porous}, left). See Table~\ref{tbl:upscaled-diffusion} for the resulting
    effective diffusion tensor.}
  \label{fig:cell-problem-1}
\end{figure}

\begin{figure}[h!]
  \centering
  \includegraphics[width=0.8\textwidth]{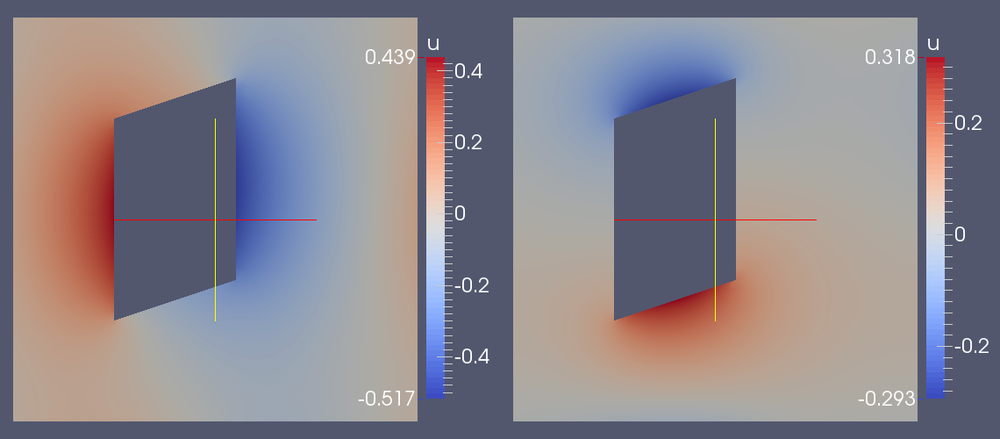}
  \caption{Solutions to the cell problems that correspond to anisotropic periodic geometry
    (Figure~\ref{fig:porous}, right). See Table~\ref{tbl:upscaled-diffusion} for the resulting
    effective diffusion tensor.}
  \label{fig:cell-problem-2}
\end{figure}

Figure  \ref{fig:cell-problem-1} and Figure \ref{fig:cell-problem-2} show the
solutions to the cell problems (\ref{eq:cell-problem}) for the
isotropic and anisotropic geometry case, respectively.  The 2D solver
for elliptic PDE with periodic boundary conditions needed for these
periodic cell problems was implemented in C++ using \texttt{deal.II} Numerics
library; see  \cite{dealII81} for details on this platform.

\begin{table}[h!]
  \centering
  \begin{tabular}{cc}
    Isotropic & Anisotropic\\\hline\\
    $\bar{\mathbb{D}}_1=\begin{bmatrix}\tensorIso\end{bmatrix}$ & $\bar{\mathbb{D}}_1=\begin{bmatrix}\tensorAniso\end{bmatrix}$
  \end{tabular}
  \caption{Examples of effective diffusion tensors corresponding to the first species (i.e. to the monomer population) for the two choices of  microstructures shown in Figure \ref{fig:porous}.}
  \label{tbl:upscaled-diffusion}
\end{table}

Controlling the cell functions allows us also to approximate the
tortuosity tensor in a direct manner, avoiding complex analytical
calculations hard to justify theoretically or experimentally; compare e.g. with Ref.
\cite{Guo}. An example in this sense is shown in Figure \ref{tortuosity}. To obtain
it, we use the relation $$\bar{\mathbb{D}}_1=d_1\phi
\bar{\mathbb{T}}^*$$ (see \cite{Bear}, e.g.) and the fact that for the
microstructures shown in Figure \ref{fig:porous} we know that the
porosity for the isotropic case is $0.75$, while the porosity for the
anisotropic case amounts to $0.85$.  We refer the reader to
\cite{Ijioma} for more numerical examples of multiscale investigations
of anisotropy effects on transport in periodically perforated media.

\begin{table}[h!]
  \centering
  \begin{tabular}{cc}
    Isotropic & Anisotropic\\\hline\\
    $\bar{\mathbb{T}}^*=\begin{bmatrix}\tensorIsotor\end{bmatrix}$ & $\bar{\mathbb{T}}^*=\begin{bmatrix}\tensorAnisotor\end{bmatrix}$
  \end{tabular}
  \caption{Examples of effective tortuosity tensors corresponding to the first species (i.e. the monomer population) for the two choices of  microstructures shown in Figure \ref{fig:porous}.}
  \label{tortuosity}
\end{table}

As soon as
the covering with microstructures lacks ergodicity and/or
stationarity, such evaluations are often replaced by efforts to
calculate accurate upper bounds on the prominent effective coefficients; see Ref. \cite{Vernescu}, for
instance, for details in this direction.

\subsection{Extensions to non-periodic microstructures}\label{extensions}

One can relax the periodicity assumption on the distribution of the
microstructures. Instead of promoting the stochastic homogenization
approach (cf. Ref. \cite{Zhikov}, e.g.) which is prohibitory expensive
from the computational point of view, we indicate two computationally
tractable cases: (1) the locally periodic arrays of microstructures
(see \cite{Boyaval_MMS,fatima2011homogenization,EJM}) and (2) the
weakly stochastic case (see \cite{LeBris3} and references cited
therein).  We will show elsewhere not only how our model formulation
and asymptotics as $\varepsilon\to 0$ translate into the frameworks of
these two non-periodic settings, but also the way the new effective
transport coefficients can be approximated numerically.


\section{Simulation studies}
\label{sec:numerics}
In this section, we study how aggregation affects deposition during the transport
of colloids in porous media. Within this frame we work with a reference parameter regime pointing out to the  {\em fast aggregation -- slow deposition regime}, that is high $\Lambda$ and low $\Biot$.

We take the model from \cite{johnson1996colloid} as the starting point
of this discussion and aim at recovering their results. We interpret
all coefficients from \cite{johnson1996colloid} in terms of our
effective coefficients obtained by the asymptotic homogenization
performed in Section \ref{sec:macroscopic-model}. As main task, we
search for new effects coming into play due to colloids aggregation.

The model for the evolution of the single mobile colloid species
$n(x,t)$ and the surface coverage of the porous matrix by the immobile
colloids $\theta(x,t)$ (that corresponds to the amount of mass
deposited) is as follows: Find the pair $(n,\theta)$ satisfying the balance equations
\begin{align}
  \label{eq:transport-mobile}
  &\partial _t n=-v_p \cdot  \nabla n + D_h\Delta n - \frac{f}{\pi a_p^2}\partial _t\theta,\\
  \label{eq:transport-immobile}
  &\partial _t \theta=\pi a_p^2 k n B(\theta),\\
  \intertext{with the switch boundary conditions}
  &n(t,0)=
    \begin{cases}
      n_0& t\in [0,t_0]\\
      0  & t>t_0
    \end{cases},\\
  &\frac{\partial n}{\partial \nu}(t,L) = 0,\\
  \intertext{and initial conditions}
  \label{eq:transport-ic-mobile}
  &n(0,x)=0,\\
  \label{eq:transport-ic-immobile}
  &\theta(0,x)=0, \ x\in [0,L].
\end{align}
Here $v_p$ is the interstitial particle velocity of the suspended
colloids, $D_h$ is the hydrodynamic particle dispersion, $a_p$ is the
particle radius, while $f$ is the specific surface area. $t_0$ is the switching off time in the boundary condition.

Given a column of cross-section surface $S$ and height $Z$ randomly
packed with spherical collector beads of radius $a_c$ and porosity (void volume
fraction) $\porosity$ typically of order of $0.4$, $f$ can be calculated
(cf. \cite{privman1991particle}) as the ratio of the total surface
area of all beads in the column to the void volume $\porosity ZS$. For
spherical beads of uniform radius, the specific surface area $f$ is
\begin{equation}
  \label{eq:specific-surface-area}
  f(\porosity):=\frac{3(1-\porosity)}{\porosity a_c}.
\end{equation}
The dynamic blocking function $B(\theta)$ arising in
(\ref{eq:transport-immobile}) accounts for the transient rate of
particle deposition.  As the colloids accumulate on the surface of the
porous matrix, they exclude a part of the surface, limiting the amount
of sites for further particle attachment.

\begin{table}[h!]
  \centering
  \begin{tabular}{p{200pt}|l}
    \hline\\
    Interstitial particle velocity                & $v_p=\frac{U}{\porosity }(2-(1-\frac{a_p}{r_0})^2)$\\
    Hydrodynamic dispersion coefficient           & $D_h=\frac{D_\infty}{\tau}+\alpha_Lv_p$\\
    Particle radius                               & $a_p=0.15\,[\mu m]$\\
    Specific surface area                         & $f=\frac{3(1-\porosity)}{\porosity a_c}$\\
    Collector grain radius                        & $a_c=0.16\,[mm]$\\
    Pore radius                                   & $r_0=(1.1969\varepsilon -0.1557)a_c$\\
    Darcy velocity                                & $U=1.02\times 10^{-4}\,[m/s]$\\
    Porosity                                      & $\porosity =0.392\,[-]$\\
    Dispersivity parameter                        & $\alpha_L=0.692\,[mm]$\\
    Kinetic rate constant                         & $k=0.25\eta U = 5\times 10^{-7}\,[m/s]$\\
    Characteristic length                         & $L=0.101\,[m]$\\
    Characteristic time                           & $t_0=5445\,[s]$\\
    Initial concentration                         & $n_0=5.58\times 10^8\,[cm^{-3}]$\\
    \\
    \hline
  \end{tabular}
  \caption{Reference parameters for simulation studies. The numerical values are taken from \cite{johnson1996colloid}.}
  \label{tab:physical-parameters}
\end{table}

We used the Finite Element Numerics toolbox DUNE
\cite{bastian2008generic} to implement a solver for the model.  We employed the
Newton method to deal with the nonlinearities in the aggregation term (counterpart of $R(\cdot)$ cf. Section \ref{aggregation})
and in the blocking function term (here denoted by $B(\cdot)$).
An implicit Euler iteration is used for time-stepping.

\begin{figure}[h!]
  \centering
  \includegraphics[width=0.8\textwidth]{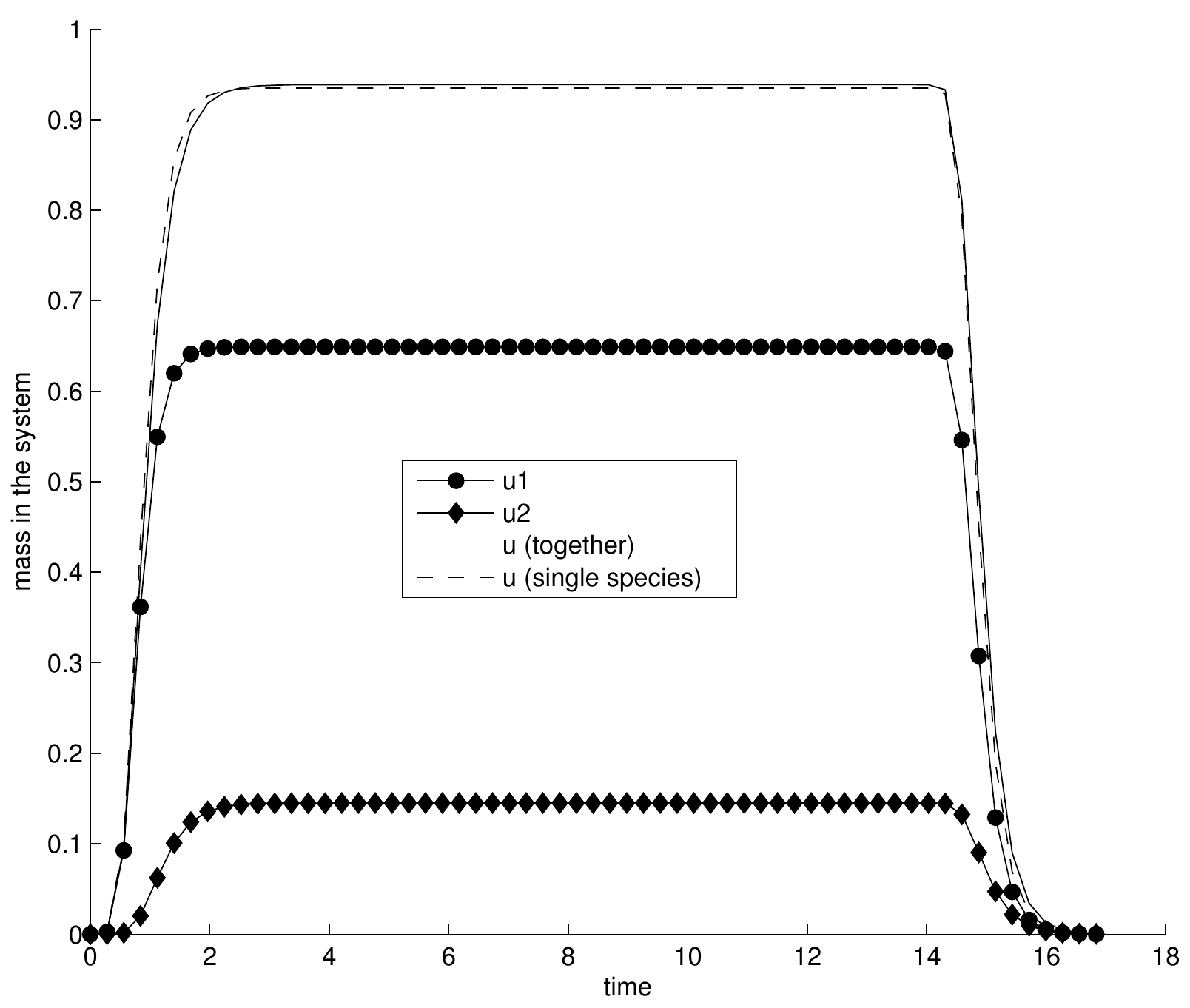}
  \caption{Simulation comparison for a single species system \emph{versus} an
    aggregating system.  The straight line is the breakthrough curve
    for the colloidal mass for the problem without aggregation.  The
    dashed line is the breakthrough curve for the colloidal mass for
    the problem with aggregation.  It is obtained by summing mass-wise
    the breakthrough curves for the monomers $u_1$ and dimers $u_2$.}
  \label{fig:simulation-comparison}
\end{figure}

The first results of our simulation with the reference parameters
indicated in Table \ref{tab:physical-parameters} are shown in
Figure~\ref{fig:simulation-comparison}.  Essentially, a single-species
system (\ref{eq:transport-mobile})-(\ref{eq:transport-ic-immobile}) is
compared to a two-species system with a square pulse going from one
side of the domain for a fixed amount of time in the first species
only.  The resulting breakthrough curves are plotted.  It is of
interest to compare the breakthrough curves for the total amount of
mass going through, no matter if it's in the form of small or large
particles. As we can observe, there is a perceptible difference
between the two curves, being the mass for the two-species case coming
in slower. This is due to larger particles having higher affinity for
deposition.

Let us focus now our attention on a specific aspect of the deposition
process, namely on the effect of the dynamic blocking functions. The
context is as follows: The rate of colloidal deposition is known to go
down as more particles attach themselves the the favorable deposition
sites of the porous matrix; see, for instance, \cite{liu1995colloid}
and references cited therein.

One of the choices for the blocking function in
(\ref{eq:transport-immobile}) corresponds to Langmuir's molecular
adsorption model \cite{langmuir1918adsorption}.  It is an affine function in
 terms of $\theta$, reaching the maximum of $1$ when the fraction of the surface
covered is zero. In other words, $B(\cdot)$ is defined as
\begin{equation}
  \label{eq:langmuir-blocking-function}
  B(\theta):=1-\beta \theta.
\end{equation}

For the simulations, we used the value  $\beta=2.9$. This corresponds to the
hard sphere jamming limit $\theta_\infty=0.345$, which is specific to
spherical collector geometry and the experimental conditions described in \cite{johnson1995dynamics}.

\begin{figure}[h!]
  \centering
  \includegraphics[width=0.8\textwidth]{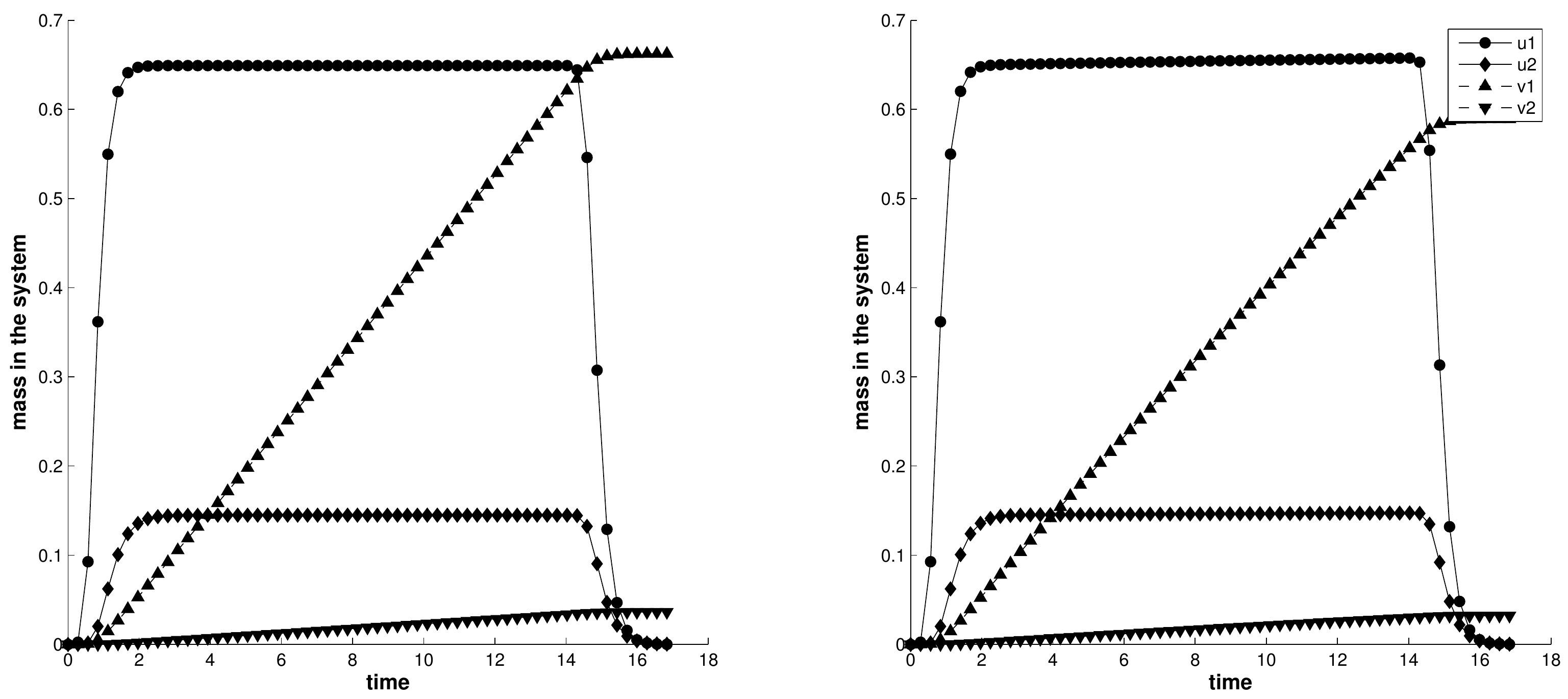}
  \caption{The effect of the Langmuirian dynamic blocking function on the
    deposition (right) \emph{versus} no blocking function (left). $u_1$ and
    $u_2$ are the breakthrough curves, while $v_1$ and $v_2$ are the
    concentrations of the deposited species.}
  \label{fig:langmuir-blocking}
\end{figure}

A simulation example of our balance equations
(\ref{eq:transport-mobile})-(\ref{eq:transport-ic-immobile}) with the
Langmuirian blocking function is shown in
Figure~\ref{fig:langmuir-blocking}.

\begin{figure}[h!]
  \centering
  \includegraphics[width=0.8\textwidth]{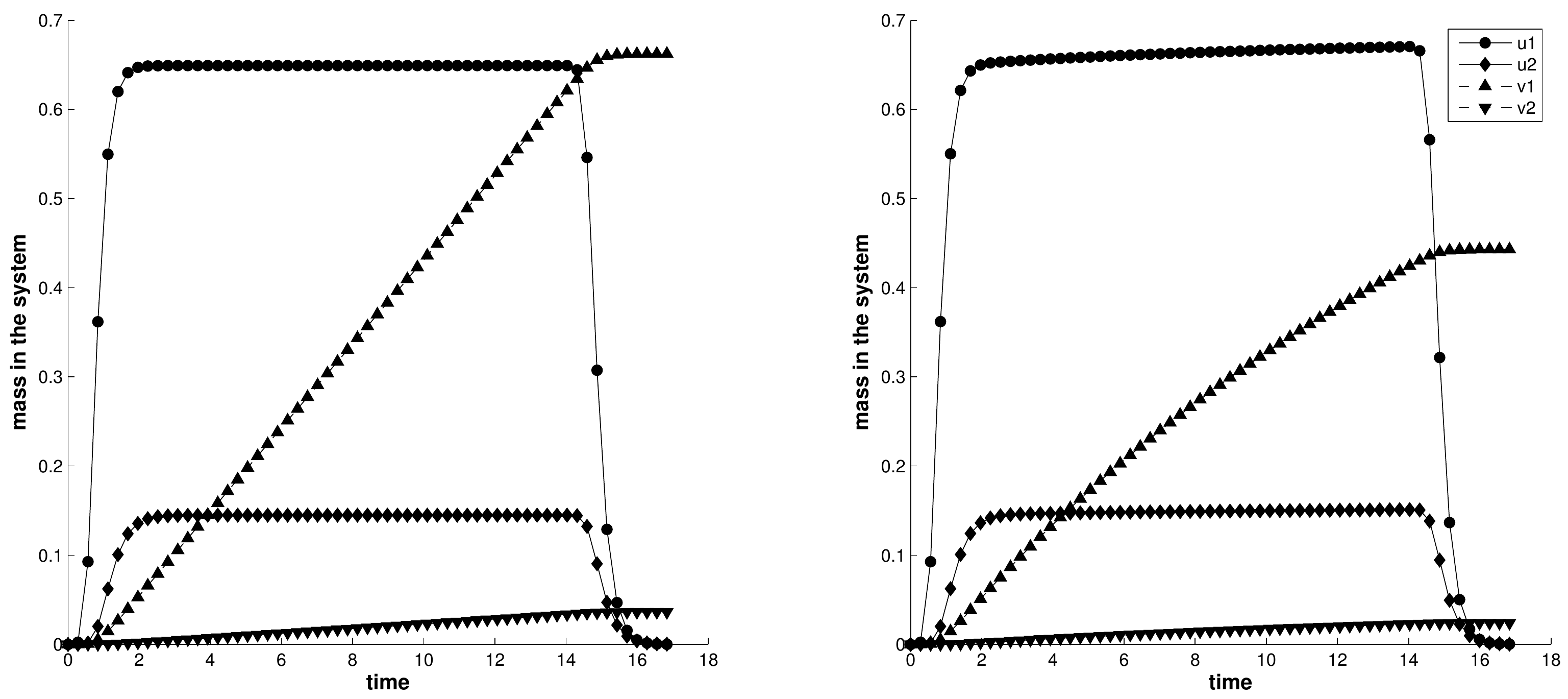}
  \caption{The effect of the RSA dynamic blocking function on the
    deposition (right) \emph{versus} no blocking function (left). $u_1$ and
    $u_2$ are the breakthrough curves, while $v_1$ and $v_2$ are the
    concentrations of the deposited species.}
  \label{fig:rsa-blocking}
\end{figure}

Another choice is the RSA dynamic blocking function as developed in
\cite{schaaf1989surface}. RSA stands for "random sequential adsorption". The RSA blocking choice is based on a third order expansion of
excluded area effects and can be used for low and moderate surface
coverage. Here $B(\theta$) is defined as:
\newcommand\thetam{\theta_\infty\beta\theta}
\begin{equation}
  \label{eq:rsa-blocking-function}
  B(\theta):=1-4\thetam+3.308(\thetam)^2+1.4069(\thetam)^3.
\end{equation}
Here, $\theta_\infty$ is the hard sphere jamming limit.  A simulation
example of the balance equations
(\ref{eq:transport-mobile})-(\ref{eq:transport-ic-immobile}) including the
RSA blocking function is shown in Figure~\ref{fig:rsa-blocking}.

\begin{figure}[h!]
  \centering
  \includegraphics[width=0.8\textwidth]{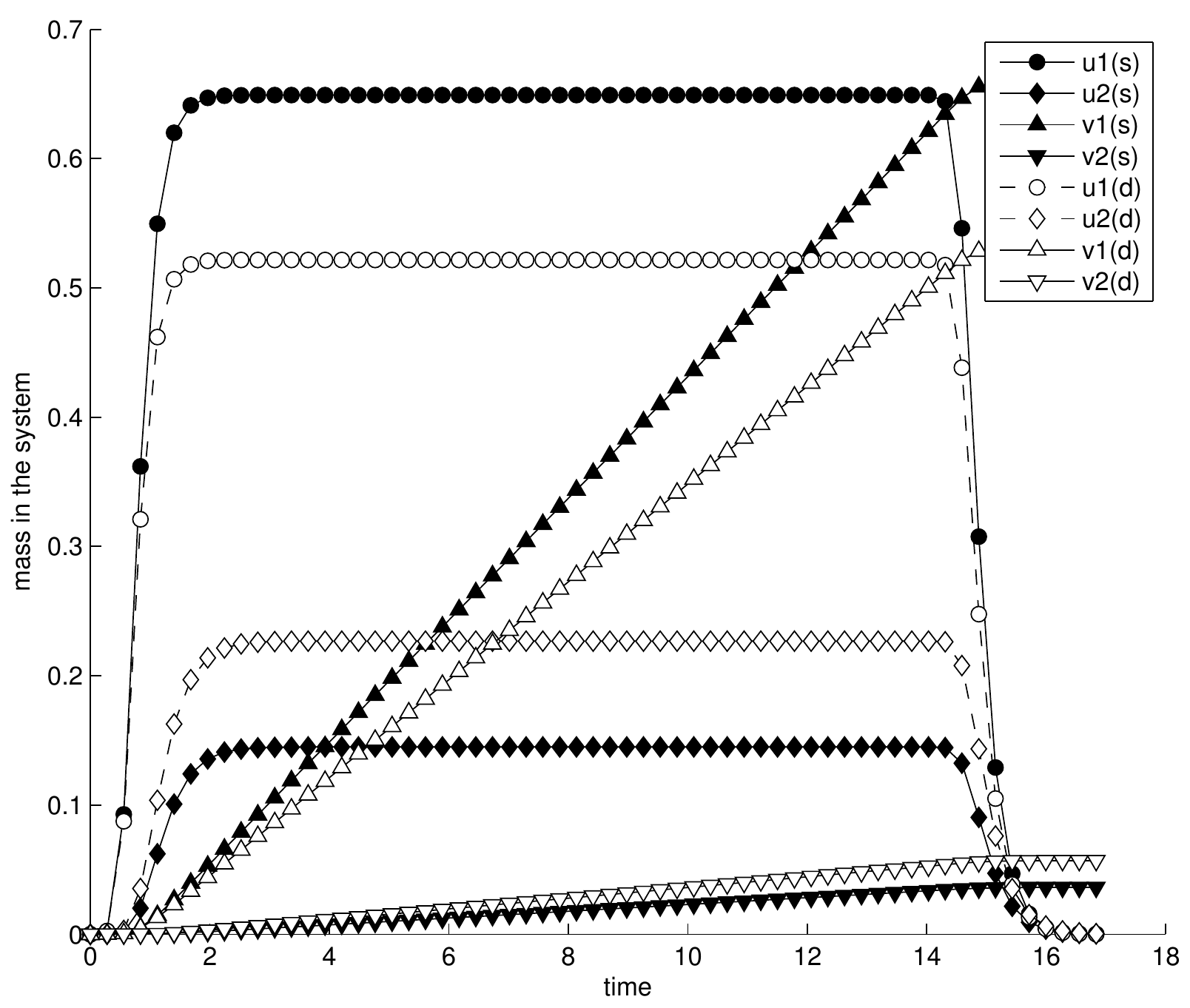}
  \caption{The effect of aggregation rates on the breakthrough
    curves. On the left, the default rate of aggregation is used, on
    the right - it's doubled. A change of aggregation rate can be
    achieved by varying the concentration of salt in the suspension,
    according to DLVO theory. Note the strong effect of aggregation on deposition. }
  \label{fig:aggregation-comparison}
\end{figure}

\section{Discussion}\label{Discussion}

This paper sheds light on  transport, aggregation/flocculation, and
deposition of colloidal particles in heterogeneous media. We succeeded
to recover basic results obtained with standard models for (single
class, single species) colloidal transport.  Furthermore, our model
includes information about the multiscale structure of the porous
medium and demonstrates new effects attributed to flocculation, such
as the occurrence of an overall decrease in the species mobility due
to a higher affinity for deposition of the large size classes of
colloidal species;  see Figure \ref{fig:aggregation-comparison} for this effect.

Extensions of this work can go in multiple directions:
\begin{itemize}
\item[(i)] Cf. \cite{liu1995colloid}, the extent of colloidal
  transport in groundwater is largely determined by the rate at which
  colloids deposit on stationary grain surfaces.  The assumption of
  stationarity can be potentially relaxed, thus aiming to incorporate
  the interplay between biofilms growth and deposition, hence
  obtaining a better understanding of the clogging/blocking of the
  pores; see e.g. \cite{SorinWRR,KnabnerZAMM}.
\item[(ii)] If repulsive forces between colloids are absent due to suitable chemical conditions,  then
  the deposition rate tends to
  increase as colloids accumulate on the grain surface (see
  Figure~\ref{fig:porous}).  Based on \cite{liu1995colloid}, this
  enhancement of deposition kinetics is attributed to the retained
  particles and is generally referred to as ripening. Active repulsive
  forces seem to lead to a decline in the deposition kinetics.  These
  effects could be investigated by our model, provided suitable
  modifications of the fluxes responsible for the transport of
  colloidal species are taken into account \cite{KnabnerHertz}.
\item[(iii)] The role of the electrolyte concentration (typically a
  salt, e.g. $KCl$) and the effect of the interplay between the
  electrostatic and van der Waals interactions on deposition kinetics
  can be studied by further developing the model. A few basic ideas on
  how to proceed in this case are collected, for instance, in
  \cite{ray2012rigorous}.
\item[(iv)] Non-periodic distributions of microstructures are relevant
  for practical applications. We leave as further work the extension
  of our solver towards the MsFEM approach, where cell problems are solved
  for each grid element, parametrized by the localized properties of
  the medium. We refer the reader to Section \ref{extensions} for
  comments in this direction.
\end{itemize}


\section*{Acknowledgments}

The authors
would like to thank Prof. Dr. Kai Uwe Totsche (Jena) and his group for very helpful
discussions on the complexity of the interactions and transport of colloids in soils.

AM and OK gratefully acknowledge financial
support by the European Union through the Initial Training Network
\emph{Fronts and Interfaces in Science and Technology} of the Seventh
Framework Programme (grant agreement number 238702).

\printbibliography
\end{document}